\documentclass[11pt]{article}
\newcommand{\documentdate}{27 V 2021}

\usepackage{a4wide,latexsym,graphicx,epstopdf,amsmath,varioref,xcolor}

\title{Adaptive Regularization Minimization Algorithms\\ with Non-Smooth
  Norms and Euclidean Curvature}

\author{
  S. Gratton%
  \thanks{Universit\'e de Toulouse, INP, IRIT, Toulouse, France. Email: serge.gratton@enseeiht.fr}
  ~and Ph. L. Toint%
  \thanks{NAXYS, University of Namur, Namur, Belgium. Email: philippe.toint@unamur.be}
}

\newcommand{\beqn}[1]{\begin{equation}\label{#1}}
\newcommand{\eeqn}{\end{equation}}
\newcommand{\req}[1]{(\ref{#1})}
\newcommand{\ms}{\;\;\;\;}
\newcommand{\tim}[1]{\;\; \mbox{#1} \;\;}

\newtheorem{theorem}{Theorem}[section]
\newtheorem{lemma}[theorem]{Lemma}
\newtheorem{corollary}[theorem]{Corollary}
\newcommand{\numsection}[1]{\section{#1}\setcounter{equation}{0}}

\renewcommand{\theequation}{\arabic{section}.\arabic{equation}}

\newcounter{algo}[section]
\renewcommand{\thealgo}{\thesection.\arabic{algo}}
\newcommand{\llem}[2]{\vspace{\baselineskip} 
\noindent\framebox[\textwidth]{\parbox{0.95\textwidth}{
\begin{lemma} \label{#1} \rm #2 \end{lemma} } } \vspace{\baselineskip} }
\newcommand{\lcor}[2]{\vspace{\baselineskip} 
\noindent\framebox[\textwidth]{\parbox{0.95\textwidth}{
\begin{corollary} \label{#1} \rm #2 \end{corollary} } } \vspace{\baselineskip} }

\DeclareMathOperator*{\argmin}{arg\,min}
\newcommand{\algo}[3]{\refstepcounter{algo}
\begin{center}\begin{figure}[htbp]
\framebox[\textwidth]{
\parbox{0.95\textwidth} {\vspace{\topsep}
{\bf Algorithm \thealgo : #2}\label{#1}\\
\vspace*{-\topsep} \mbox{ }\\
{#3} \vspace{\topsep} }}
\end{figure}\end{center}}
\newcommand{\bpr}{{\bf Proof.} \hspace{1.5mm}}
\newcommand{\epr}{\hfill $\Box$ \vspace*{1em}}
\newcommand{\proof}[1]{
\begin{list}{}{
\setlength{\topsep}{0.0pt}
\setlength{\partopsep}{0.0pt}
\setlength{\leftmargin}{0.025\textwidth}
\setlength{\rightmargin}{0.5\leftmargin}
\setlength{\labelwidth}{0.5\leftmargin}
\setlength{\labelsep}{0.25\leftmargin}}
\item \bpr #1 \epr \noindent
\end{list}}
\newcommand{\lthm}[2]{\vspace{\baselineskip} 
\noindent\framebox[\textwidth]{\parbox{0.95\textwidth}{
\begin{theorem} \label{#1} \rm #2 \end{theorem} } } \vspace{\baselineskip} }
\newcommand{\pd}[1]{\langle #1 \rangle}
\newcommand{\ii}[1]{\{ 1, \ldots, #1 \}}

\newcommand{\calN}{{\cal N}} 
\newcommand{\calO}{{\cal O}} 
\newcommand{\calS}{{\cal S}}

\renewcommand{\Re}{\hbox{I\hskip -2pt R}}
\newcommand{\bigfrac}[2]{\frac{\displaystyle #1}{\displaystyle #2}}
\newcommand{\smallRe}{\hbox{\footnotesize I\hskip -2pt R}}
\newcommand{\sfrac}[2]{{\scriptstyle \frac{#1}{#2}}}
\newcommand{\half}{\sfrac{1}{2}}
\newcommand{\third}{\sfrac{1}{3}}

\newcommand{\sixth}{\sfrac{1}{6}}
\newcommand{\eqdef}{\stackrel{\rm def}{=}}

\newcommand{\al}[1]{{\footnotesize{\sf #1}}}
\newcommand{\tal}[1]{{\normalsize {\sf #1}}}
\newcommand{\nr}[1]{\|#1\|_r}

\date{\documentdate}

\begin{document}

\maketitle

\begin{abstract}
A regularization algorithm (\al{AR$1p$GN}) for unconstrained nonlinear minimization is
considered, which uses a model consisting of a Taylor expansion of
arbitrary degree and regularization term involving a possibly
non-smooth norm. It is shown that the non-smoothness of the
norm does not affect the $\calO(\epsilon_1^{-(p+1)/p})$ upper bound on
evaluation complexity for finding first-order $\epsilon_1$-approximate
minimizers using $p$ derivatives, and that this result does not hinge on the
equivalence of norms in $\Re^n$.  It is also shown that, if $p=2$,
the bound of $\calO(\epsilon_2^{-3})$ evaluations for finding second-order
$\epsilon_2$-approximate minimizers still holds for a variant
of \al{AR$1p$GN} named \al{AR2GN}, despite the possibly
non-smooth nature of the regularization term.  Moreover, the adaptation of the
existing theory for handling the non-smoothness results in an interesting
modification of the subproblem termination rules, leading to an even more
compact complexity analysis.
In particular, it is shown when the Newton's step is acceptable for an adaptive
regularization method. The approximate minimization of quadratic polynomials regularized
with non-smooth norms is then discussed, and a new approximate second-order necessary
optimality condition is derived for this case.  An specialized algorithm is then
proposed to enforce the first- and  second-order conditions that are strong enough
to ensure the existence of a suitable step in \al{AR$1p$GN} (when $p=2$) and
in \al{AR2GN}, and its iteration complexity is analyzed.
\end{abstract}

{\small

\vspace*{3mm}
  \textbf{Note:} This paper is a close variant of \cite{GratToin21} where the theory presented
  there is adapted to the use of Euclidean curvature for second-order
  optimality conditions, instead of the (potentially hard to compute)
  most negative curvature in a non-smooth norm. We have chosen to keep the complete
  presentation to preserve its self-contained character and because the
  necessary modifications to \cite{GratToin21} are scattered throughout the text.
}
\numsection{Introduction}

This paper is concerned with the derivation of upper bounds on the evaluation
complexity of adaptive regularization algorithms for the solution of the
unconstrained nonconvex optimization problem \beqn{problem} \min_{x \in
  \smallRe^n} f(x).  \eeqn This research area has been remarkably active in
recent years (see, for instance,
\cite{Grie81,NestPoly06,CartGoulToin09a,CartGoulToin11d,CartGoulToin12a,%
  BianLiuzMoriScia15,GrapYuanYuan15a,BianScia16,BirgGardMartSantToin17,Mart17,%
  GratRoyeVice17,BergDiouGrat17,BellGuriMoriToin19,CartGoulToin20b}).
Adaptive regularization algorithms, the class of methods considered here,
compute steps from one iterate to the next by building and (often
approximately) minimizing a model consisting of a truncated Taylor expansion
of $f$, which is then ``regularized'' by adding a suitable power of the norm
of the putative step. Several authors have considered various smooth norms for
this regularization term
\cite{NestPoly06,CartGoulToin11d,Duss15,BirgGardMartSantToin17,ChenToinWang19,ChenToin20},
showing that, under suitable assumptions, the resulting adaptive
regularization method must find a first-order $\epsilon_1$-approximate
minimizer for problem \req{problem} (that is an iterate $x_k$ such
$\|\nabla_x^1f(x_k)\| \leq \epsilon_1$) in at most
$\calO(\epsilon_1^{-(p+1)/p})$ evaluations of the objective function and its
derivatives. In addition, second-order variants of this algorithm are bound to
find a second-order $\epsilon_2$-approximate minimizer (that is an iterate
$x_k$ such the smallest eigenvalue of $\nabla_x^2f(x_k)$ exceeds
$-\epsilon_2$) in at most $\calO(\epsilon_2^{-(p+1)/{(p-1)}})$ such
evaluations. The detailed algorithms considered in these contributions all
depend on the central tenet that the regularized model (whose approximate
minimization yields the step from one iterate to the next) is smooth, and thus
that this approximate minimization can be carried out using algorithms for
smooth functions and can be terminated using approximate optimality conditions
for smooth problems. We show in this paper that the same evaluation complexity
bounds still holds for first-order approximate minimizers in the case where
non-smooth norms (such as $\ell_1$ or $\ell_\infty$) are considered,
provided the algorithm is suitably modified. We also show that, when $p=2$,
the evaluation complexity bound in $\calO(\epsilon_2^{-3})$ is also maintained
in the same context for a variant of the algorithm. Unsurprisingly, both
results require redefining the termination conditions for model
minimization. As it turns out, the resulting modifications of the standard
adaptive regularization method are extremely simple and their use in the
complexity theory results in a remarlably compact formulation.

One may argue that, since all norms are equivalent in finite dimensional
spaces, the stated complexity bound can be derived for any norm from known
results in Euclidean norm \cite{BirgGardMartSantToin17,CartGoulToin20}.  While
this is true if one focuses on the order in $\epsilon_1$ and $\epsilon_2$
only, this ignores the influence of the norm equivalence constants, whose size
can be significant when $n$, the dimension of the problem, grows.  For
instance the equivalence constant between the Euclidean and infinity norm is
proportional to the square root of the problem's dimension.  Thus obtaining a
given accuracy on the gradient norm in the infinity norm by simply applying
the norm equivalence principle may require $n^{(p+1)/2p}$ times more
evaluations of the objective function and its derivatives than in the
Euclidean one. The approach presented here attempts to avoid this
potentially problematic increase in cost.

Of course, for the new algorithms to be practical, one needs to show
that the model minimization subproblems are solvable by implementable
methods. Focusing again on the case where $p=2$ and the model is a regularized
quadratic, we derive a specialized second-order necessary optimality
condition for the approximate minimization of such non-smooth functions.  We then
propose a new algorithm which is able to achieve first- and second-order
approximate optimality for this problem and evaluate its iteration complexity.
We finally discuss relaxed variants of the new algorithms that are sufficient for
solving the subproblems of interest in algorithms for general functions, as
well as their iteration complexity. 

Our exposition is organized as follows. We present the problem and the
first-order algorithm in Section~\ref{algo-s} and derive its evaluation
complexity theory in Section~\ref{complexity-s}. Section~\ref{complexity2-s}
discusses the new approximate second-order necessary condition for global
minimizers and establishes the upper bound on evaluation complexity for an
adapted variant of the algorithm.  A method for approximately minimizing
regularized quadratics (enough for solving the subproblems 
arising in Sections~\ref{algo-s} and \ref{complexity2-s})
is then presented and analyzed in Section~\ref{RQMIN-s}.  Finally, a brief
conclusion is stated in Section~\ref{conclusion-s}.

\numsection{An first-order adaptive regularization in general norms}
\label{algo-s}

\noindent
We consider the unconstrained nonlinear optimization problem \req{problem}
where $f$ is a (potentially nonconvex) $p$ times continuoulsy differentiable
function from $\Re^n$ to $\Re$, for some integer $p\geq 1$.  We define
\[
T_{f,p}(x,s)
\eqdef f(x) + \sum_{\ell=1}^p \frac{1}{\ell !} \nabla_x^\ell f(x)[s]^\ell,
\]
the $p$-th order Taylor expansion of $f$ at $x$, where the notation
$\nabla_x^\ell f(x)[s]^\ell$ denotes the symmetric $\ell$-dimensional tensor
$\nabla_x^\ell f(x)$ applied on $\ell$ copies of the vector $s$.

As outlined in the introduction, adaptive regularization methods
are iterative schemes that compute a step form an
iterate $x_k$ by constructing a regularized model $m_k(s)$ of $f(x_k+s)$ as 
\beqn{model}
m_k(s) \eqdef T_{f,p}(x_k,s) + \frac{\sigma_k}{(p+1)!} \nr{s}^{p+1},
\eeqn
where the $p$-th order Taylor series is ``regularized'' by adding the term
$ \sigma_k \nr{s}^{p+1}/(p+1)!$ ($\sigma_k$ is known as the ``regularization
parameter'') and where we allow $\nr{\cdot}$ to be a general possibly
non-smooth norm. This implies that $\nr{\cdot}$ is convex and Lipschitz
continuous with global Lipschitz constant equal to one. Given the
$\nr{\cdot}$ norm and defining 
\beqn{Tnorm}
\| S_j \|_{r,j} = \max_{\nr{s} = 1} |S_j[s]^j|
\eeqn
to be the norm of the $j$-dimensional symmetric (for $j>1$) tensor $S_j$  induced by 
$\nr{\cdot}$, we are now interested in finding, for some prespecified accuracy
requirement $\epsilon_1 \in (0,1]$, an
$\epsilon_1$-approximate first-order critical point, that is a point
$x_{\epsilon_1}$ such that 
$
\|\nabla_x^1 f(x_\epsilon)\|_{r,1} \leq \epsilon_1.
$
Note that, because of \req{Tnorm},
$\|\cdot\|_{r,1}$ is the dual norm of $\nr{\cdot}$.)

The ``regularization term'' in \req{model} guarantees that $m_k(s)$ is bounded
below and thus makes the procedure of finding a step $s_k$ by (approximately)
minimizing $m_k(s)$ well-defined. However, at variance with the usual setting
for adaptive regularization methods, \emph{the model $m_k(s)$  may no longer be smooth}.
Once the step is computed, the value of the objective function at the
trial point $x_k+s_k$ is then computed.  If the decrease in $f$ from $x_k$ to
$x_k+s_k$ is comparable to that predicted by the second-order Taylor series,
the trial point is accepted as the new iterate and the regularization
parameter is (possibly) reduced. If this is not the case, the trial point is
rejected and the regularization parameter increased.  The resulting algorithm
is formally stated as the \al{AR$1p$GN} algorithm \vpageref{AR1pGN}.

\algo{AR1pGN}{First-Order Adaptive Regularization with General Norm (\tal{AR$1p$GN})}{
  \begin{description}
  \item[Step 0: Initialization. ] An initial point $x_0\in \Re^n$, a regularization
    parameter $\sigma_0$ and a desired final gradient accuracy
    $\epsilon_1 \in (0,1]$ are given. The constants
    $\eta_1$, $\eta_2$, $\gamma_1$, $\gamma_2$, $\gamma_3$, $\theta_1$ and $\sigma_{\min}$
    are also given such that
    \beqn{eta-gamma}
    \sigma_{\min} \in (0, \sigma_0], \ms 
      0 < \eta_1 \leq \eta_2 < 1, \ms \theta_1 > 1
      \tim{ and } 0< \gamma_1 < 1 < \gamma_2 < \gamma_3.
    \eeqn
    Compute $f(x_0)$ and set $k=0$.
  \item[Step 1: Check for termination. ] Terminate with $x_{\epsilon_1} = x_k$ if
    \beqn{term}
    \|\nabla_x^1 f(x_k)\|_{r,1} \leq \epsilon_1.
    \eeqn
  \item[Step 2: Step calculation. ] Compute a step $s_k$  which
    sufficiently reduces the model $m_k$ in the sense that
    \beqn{descent}
    m_k(s_k) \leq m_k(0)
    \eeqn
    and
    \beqn{term-model}
    \|\nabla_s^1 T_{f,p}(x_k,s_k)\|_{r,1} \leq \theta_1 \frac{\sigma_k}{p!}\nr{s}^p.
    \eeqn
  \item[Step 3: Acceptance of the trial point. ]
    Compute $f(x_k+s_k)$ and define 
    \beqn{rhok-def}
    \rho_k = \frac{f(x_k) - f(x_k+s_k)}{T_{f,p}(x_k,0)-T_{f,p}(x_k,s_k)}.
    \eeqn
    If $\rho_k \geq \eta_1$, then define
    $x_{k+1} = x_k + s_k$; otherwise define $x_{k+1} = x_k$.
  \item[Step 4: Regularization parameter update. ]
   Set
   \beqn{sigma-update}
   \sigma_{k+1} \in \left\{ \begin{array}{ll}
   {}[\max(\sigma_{\min}, \gamma_1\sigma_k), \sigma_k ]  & \tim{if} \rho_k \geq \eta_2, \\
   {}[\sigma_k, \gamma_2 \sigma_k ]          &\tim{if} \rho_k \in [\eta_1,\eta_2),\\
   {}[\gamma_2 \sigma_k, \gamma_3 \sigma_k ] & \tim{if} \rho_k < \eta_1.
   \end{array} \right.
   \eeqn
   Increment $k$ by one and go to Step~1.
  \end{description}
}

\noindent
While the \al{AR$1p$GN} algorithm follows the main lines of existing adaptive
regularization methods (see \cite{CartGoulToin11d,BirgGardMartSantToin17} for
example), we immediately note that the test \req{term-model} differs from the
test $\| \nabla_s^1 m_s(s_k) \|_2 \leq \theta_1 \|s_k\|_2^p$ which is used so
far in the literature.  Indeed, our framework no longer guarantees that
$\nabla_s^1 m_s(s)$ exists, due to the possible lack of smoothness of the
regularization term. Note however that, if $\nr{\cdot}$ is differentiable
everywhere except at the origin, then
\[
\nabla_s^1 T_{f,p}(x_k,s_k) + \frac{\sigma_k}{p!}\nr{s}^{p}\, \nabla_s^1\nr{s}
=0
\]
at a nonzero first-order point of $m_k(s)$, and \req{term-model} holds at such a point
since $\nr{s}$ is Lipschitz continuous with unit Lipschitz constant, and thus
$\|\nabla_s^1\nr{s}\|\leq 1$. The condition \req{term-model} is
therefore weaker than a more usual condition of the form $\|\nabla_s^1
m_k(s_k)\|\leq \theta_1 \|s_k\|^p$.

Remarkably, and at variance with other adaptive regularization methods, the
\al{AR$1p$GN} algorithm allows the Newton step $s_k =
-\nabla_x^2f(x_k)^{-1}\nabla_x^1 f(x_k)$ when $p=2$ and the  Hessian $\nabla_x^2f(x_k)$
is positive definite, \emph{provided the regularized model has not increased, that
is provided \req{descent} holds}. Indeed this step automatically ensures
\req{term-model} since then $\nabla_s^1 T_{f,2}(x_k,s_k)=0$. The condition
\req{descent} however avoids situations where the model decrease $m_k(0)-m_k(s_k)$
is tiny but $\nr{s_k}$ is large, which is exactly what happens in the example
of \cite{CartGoulToin18a} showing convergence of Newton's method to a first-order
$\epsilon_1$-approximate minimizer in $\calO(\epsilon_1^{-2})$ evaluations.

We also note that we could use an iteration-dependent $\theta_{1,k}$ in
\req{term-model}, provided it is bounded below by one and strictly bounded above by a
constant. We have ignored this possibility for the sake of simplicity.

Having modified the requirements on the step, we now need to verify that the
new conditions \req{descent} and \req{term-model} are compatible. We start by
deriving an expression for the subdifferential $\partial( \nr{\cdot}^{p+1})(s)$.

\llem{power-subdiff}{
  We have that 
  \beqn{subdiff-norm}
  \partial( \nr{\cdot})(s) 
  = \{ v \in \Re^n \mid v^Ts = \nr{s} \tim{ and } \|v\|_{r,1} = 1 \}
  \eeqn
  and
  \beqn{subdiff-power-subdiff}
  \partial_C( \nr{\cdot}^{p+1})(s)
  = \partial( \nr{\cdot}^{p+1})(s) = (p+1) \nr{s}^p \, \partial (\nr{\cdot})(s)
  \eeqn
  where $\partial_C$ denotes the Clarke subdifferential.
}

\proof{
  The identity \req{subdiff-norm} is standard (see
  \cite[Example~3.1]{HirrLema96} for instance). By composition of the norm
  with the increasing convex differentiable function $\phi(t)=t^{p+1}$
  (on $\Re^+$), we obtain from \cite[Theorem~4.3.1]{HirrLema96} that
  \[
  \partial( \nr{\cdot}^{p+1})(s) = \{ \alpha s\in \Re^ns \mid  \alpha \in
  \phi'(\nr{s}) \tim{and} s \in \partial ( \nr{\cdot})(s) \} .
  \]
  which is the second equality in \req{subdiff-power-subdiff}.
  Since $\nr{\cdot}^{p+1}$ is also Lipschitz continuous, it is Clarke
  regular and thus the Clarke subdifferential and the standard one coincide
  (see \cite[Proposition~4.3]{ClarLedySterWole98}), giving the first equality
  in \req{subdiff-power-subdiff}.
} 

\noindent
This allows us to derive the following characterization of a minimizer of $m_k$.

\llem{noc-1st-order}{Let $s_k^*$ be a local minimizer of $m_k$.  Then
  \beqn{noc-1st}
  \|\nabla_s^1T_{f,p}(x_k,s_k^*)\|_{r,1} = \half \sigma \nr{s_k^*}^p.
  \eeqn
}

\proof{
   Since $m_k$ is Lipschitz continuous, the Clarke criticality of $s_k^*$
   implies that
   \beqn{clarke-critical}
   0 \in \partial_C m_k(s_k^*)
   = \left \{ \nabla_x^1 T_{f,p}(x_k,s_k^*) \right \}+\frac{\sigma_k}{(p+1)!}\partial_C(\nr{s_k^*}^{p+1}),
   \eeqn
   where we have used the property of the Clarke subdifferential of the sum of
   two locally Lipschitz functions \cite[Exercice~1.4]{ClarLedySterWole98} and
   the fact that, since $T_{f,p}(x_k,s)$ is continuously differentiable as a
   function of $s$, $\partial_CT_{f,p}(x,.)(s) = \{\nabla_s^1
   T_{f,p}(x_k,s)\}$. Using now \req{subdiff-power-subdiff}, we deduce from
   this identity and \req{clarke-critical} that there exists a vector
   $\xi \in  \partial (\nr{\cdot})(s_k^*)$ such that 
   \beqn{noc1-1}
   \nabla_s^1 T_{f,p}(x_k,s_k^*) = - \frac{\sigma_k}{p!} \nr{s_k^*}^p \  \xi.
   \eeqn
   Moreover, \req{subdiff-norm} implies that $\|\xi\|_{r,1} = 1$. Taking norms in
   this relation gives \req{noc-1st}.
} 

\noindent
The (scalar) necessary condition \req{noc-1st} is clearly weaker that the
(vector) identity \req{clarke-critical}, but is nevertheless sufficient to
derive the following crucial result.

\lcor{step-exists}{
  A step satisfying both \req{descent}
  and \req{term-model} always exists.
}

\proof{
   From
   \[
   m_k(s)
  \geq \frac{\sigma_k}{(p+1)!}\nr{s}^{p+1}
      -|f(x)|-\sum_{\ell=1}^p\frac{1}{\ell!}\|\nabla_x^\ell f(x_k)\|_{r,\ell}\nr{s} ^\ell
   \]
   we obtain  $\lim_{\nr{s} \rightarrow + \infty} m_k(s) = +\infty$ which,
   together with the continuity of $m_k(s)$, implies that $m_k$
   admits at least one minimizer $s_k^*$ over $\Re^n$, satisfying
   $m_k(s_k^*)\leq m_k(0)$. Applying Lemma~\ref{noc-1st-order} then gives that
   \req{term-model} holds at $s_k^*$ for any $\theta_1 \geq 1$.
} 

\noindent
An important comment is in order at this point. Because the Clarke
subdifferential of the norm is not necessarily continuous in our context, it
may seem at first sight that obtaining a step satisfying the conditions \req{descent}
and \req{term-model} may require the computation of an exact minimizer $s_k^*$
of the model, which is potentially costly.  Fortunately, this fear is unfounded
because both the left- and the right-hand sides of \req{term-model} are
continuous functions of $s$ and the inequality therefore also holds in a
neighbourhood of $s_k^*$ provided $\theta_1>1$. Any convergent minimization
algorithm (such as those proposed, for instance, in
\cite{Flet82b,Nest13,LewiWrig16,CartGoulToin11a,CartGoulToin20a,GratSimoToin20} or, more
generally, in \cite{Knos17}, or in Section~\ref{RQMIN-s})
applied on the model is therefore bound to produce a suitable step $s_k$ in a
finite number of iterations.

Following well-established practice, we now define
\[
\calS \eqdef \{ k \geq 0 \mid x_{k+1} = x_k+s_k \} = \{ k \geq 0 \mid \rho_k \geq \eta_1 \}
\tim{and}
\calS_k \eqdef \calS \cap \ii{k},
\]
the set of indeces of ``successful iterations'', and
the set of indeces of successful iterations up to iteration $k$, respectively.
We also recall a well-known result bounding the total number of iterations of
an adpative regularization method in terms of the number of successful ones.

\llem{SvsU}{
Suppose that the \al{AR$1p$GN} algorithm is used and that $\sigma_k \leq
\sigma_{\max}$ for some $\sigma_{\max} >0$. Then
\beqn{unsucc-neg}
k \leq |\calS_k| \left(1+\frac{|\log\gamma_1|}{\log\gamma_2}\right)+
\frac{1}{\log\gamma_2}\log\left(\frac{\sigma_{\max}}{\sigma_0}\right).
\eeqn
}

\proof{See \cite[Theorem~2.4]{BirgGardMartSantToin17}.}

\numsection{Evaluation complexity for the \tal{AR$1p$GN} algorithm}\label{complexity-s}

\noindent
Before discussing our analysis of evaluation complexity, we first formalize
our assumptions on problem \req{problem}. 

\noindent
\fbox{\bf AS.1} $f$ is $p$ times differentiable and its $p$-th
derivative $\nabla_x^p f(x)$ is is globally Lipschitz continuous in the
$\nr{\cdot}$ and $\|\cdot\|_2$ norms, that is there exists $L_{r,p}, L_{2,p}\geq 0$ such that 
\beqn{rLip}
\|\nabla_x^p f(x)-\nabla_x^p f(y)\|_{r,p} \leq L_{r,p}\nr{x-y} \tim{ for all } x,y \in \Re^n,
\eeqn
where the $\nr{\cdot}$ norm in the left-hand side is defined by \req{Tnorm}, and
\beqn{2Lip}
\|\nabla_x^p f(x)-\nabla_x^p f(y)\|_2 \leq L_{2,p}\|x-y\|_2 \tim{ for all } x,y \in \Re^n.
\eeqn

\noindent
\fbox{\bf AS.2} There exists a constant $f_{\rm low}$ such that
$f(x) \geq f_{\rm low}$ for all $x\in\Re^n$.

\noindent
Assumption AS.2 ensures that problem \req{problem} is well-defined.
Assumption AS.1 recasts the usual context for the analysis of complexity of
adaptive regularization methods in the context of the general norms, 
Note that assuming both \req{rLip} and \req{2Lip} in AS.1 is important to avoid large
dimension-dependent ``norm-equivalence'' constants in our
final evaluation complexity bounds. Thus, for the bounds to be meaningful, we
implicitly assume that $L_{2,p}$ and $L_{r,p}$ do not vary too much in size.
AS.1 yields the well-known Lipschitz error
bounds.

\llem{lipschitz}{
  Suppose that AS.1 holds and that $k\in \calS$. Then
  \beqn{Lip-f}
  |f(x_{k+1})-T_{f,p}(x_k,s_k))| \leq \frac{L_{r,p}}{(p+1)!}\nr{s_k}^{p+1},
  \eeqn
  \beqn{Lip-g}
  \|\nabla_x^1 f(x_{k+1}) - \nabla_s^1 T_{f,p}(x_k,s_k)\|_{r,1}  \leq \frac{L_{r,p}}{p!} \nr{s_k}^p
  \eeqn
  and
  \beqn{Lip-H}
  \|\nabla_x^2 f(x_{k+1}) - \nabla_s^2 T_{f,p}(x_k,s_k)\|_2  \leq \frac{L_{2,p}}{(p-1)!} \|s_k\|_2^{p-1}.
  \eeqn
}

\proof{
  The proof of \req{Lip-f} and \req{Lip-g} is a direct extension
  of \cite[Lemma~2.1]{CartGoulToin20b}
  with $\beta=1$ that now uses \req{rLip} instead of \req{2Lip} and exploits
  \req{Tnorm}. It is given in the appendix of \cite{GratToin21}. The
  inequality \req{Lip-H} immediately results from \req{2Lip} and 
  \cite[Lemma~2.1]{CartGoulToin20b}.
}

From now on, the analysis in this section follows that presented in
\cite{BirgGardMartSantToin17} quite closely. We first state a simple lower
bound on the decrease of the Taylor expansion.
    
\llem{model-decrease}{
  \beqn{Tdecr}
  \Delta T_{f,p}(x_k,s_k) \eqdef T_{f,p}(x_k,0)-T_{f,p}(x_k,s)
  \geq \frac{\sigma_k}{(p+1)!}\nr{s_k}^{p+1}.
  \eeqn
}

\proof{Direct from \req{descent} and \req{model}.}

\noindent
We next derive an upper bound on the regularization parameter.

\llem{sigma-max}{
  Suppose that AS.1 holds.  Then, for all $k\geq 0$,
  \beqn{sigma-upper}
  \sigma_k\leq\sigma_{\max} \eqdef \gamma_3\max\left[\sigma_0,\frac{L_{r,p}}{(1-\eta_2)}\right].
  \eeqn
}

\proof{See \cite[Lemma~2.2]{BirgGardMartSantToin17}.
Using \req{rhok-def}, \req{Lip-f},
and \req{Tdecr}, we obtain that
\[
|\rho_k-1|
\leq  \bigfrac{(p+1)!|f(x_k+s_k)-T_{f,p}(x_k,s_k)|}{\sigma_k \nr{s_k}^{p+1}}
\leq \bigfrac{L_{r,p}}{\sigma_k}.
\]
Thus, if $\sigma_k \geq L_{r,p}/(1-\eta_2)$, then $\rho_k \geq \eta_2$, iteration
$k$ is successful and \req{sigma-update} implies that
$\sigma_{k+1}\leq \sigma_k$. The mechanism of the algorithm then guarantees that
\req{sigma-upper} holds.
} 

\noindent
The next lemma remains in the spirit of
\cite[Lemma~2.3]{BirgGardMartSantToin17}, but now takes the new condition
\req{term-model} into account, avoiding any reference to the model's
derivative and resulting in a simpler proof.

\llem{useful}{
  Suppose that AS.1 holds and that $k\in\calS$ before termination.  Then
  \beqn{crucial}
  \nr{s_k}^p \geq  \frac{p!}{L_{r,p}+ \theta_1\sigma_{\max}}\,\epsilon_1.
  \eeqn
}

\proof{
  Successively using the fact that termination does not occur at iteration
  $k$, the triangle inequality, \req{Lip-g} for $j=1$, condition \req{term-model} and
  \req{sigma-upper}, we deduce that 
  \[
  \begin{array}{lcl}
    \epsilon_1 & < & \|\nabla_x^1f(x_{k+1})\|_{r,1}\\*[2ex]
    & \leq & \|\nabla_x^1f(x_{k+1})-\nabla_x^1 T_{f,p}(x_k,s_k)\|_{r,1} + \|\nabla_x^1 T_{f,p}(x_k,s_k)\|_{r,1}\\*[2ex]
    & \leq & \bigfrac{L_{r,p}}{p!}\nr{s_k}^p + \theta_1 \bigfrac{\sigma_k}{p!} \nr{s_k}^p.
  \end{array}
  \]
  This in turn directly implies \req{crucial}.
} 

\noindent
We may now resort to the classical ``telescoping sum'' argument to obtain the
desired complexity result.

\lthm{complexity}{
  Suppose that AS.1--AS.2 hold.   Then the
  \al{AR$1p$GN} algorithm requires at most
  \[
  \frac{(p+1)!}{\eta_1\sigma_{\min}}\left(\frac{L_{r,p}+\theta_1\sigma_{\max}}{p!}\right)^{\frac{p+1}{p}}
  \frac{f(x_0)-f_{\rm low}}{\epsilon_1^{\frac{p+1}{p}}}
  \]
  successful iterations and evaluations of $\{\nabla_x^i f\}_{i=1,2}$ and at most
  \[
  \frac{(p+1)!}{\eta_1\sigma_{\min}}\left(\frac{L_{r,p}+\theta_1\sigma_{\max}}{p!}\right)^{\frac{p+1}{p}}
    \frac{f(x_0)-f_{\rm low}}{\epsilon_1^{\frac{p+1}{p}}}
  \left(1+\frac{|\log\gamma_1|}{\log\gamma_2}\right)+
  \frac{1}{\log\gamma_2}\log\left(\frac{\sigma_{\max}}{\sigma_0}\right)
  \]
  evaluations of $f$ to produce a vector $x_{\epsilon_1}\in \Re^n$ such that
  $\|\nabla_x^1 f(x_{\epsilon_1})\|_{r,1}\leq \epsilon_1$.
}

\proof{
  Let $k$ be the index of an iteration before termination.  Then, using AS.2,
  the definition of successful iterations, \req{Tdecr} and \req{crucial},
  \[
  |\calS_k|
  \leq \frac{(p+1)!}{\eta_1\sigma_{\min} }
  \left(\frac{L_{r,p}+\theta_1\sigma_{\max}}{p!}\right)^{-\frac{p+1}{p}}
  \frac{f(x_0)-f_{\rm low}}{\epsilon_1^{\frac{p+1}{p}}}
  \]
  for any $k$ before termination, and the first conclusion follows since the
  derivatives are only evaluated once per successful iteration.  Applying now
  Lemma~\ref{SvsU} gives the second conclusion.
} 

\numsection{Approximate second-order minimizers for $p=2$}\label{complexity2-s}

We now turn the second-order case and from now on, limit our
analysis to the case where $p=2$. We are thus interested in finding
an approximate second-order minimizer, that is,
for our present purposes,  an iterate $x_k$ such that
\beqn{lar-def}
\lambda_{\min}[\nabla_x^2 f(x_k)] \ge - \epsilon_2,
\eeqn
In what follows, we assume that, for a symmetric $H$, we can compute
$\lambda_{\min}[H] $ and a an associated vector $u$ such that $\|u\|_2=1$.
This is at variance with the approach taken in \cite{GratToin21}, where
\req{lar-def} is replaced by the condition that
\[
\min_{\nr{v}=1}\nabla_x^2 f(x_k)[v]^2 \ge - \epsilon_2,
\]
a much more difficult task for which approximations are often necessary.

We now establish a second-order necessary condition for a global minimizer
of a regularized quadratic $m$. As a first step, we derive a lower
bound on the model decrease that can be obtained along a direction of
sufficient negative curvature.

\llem{2nd-order-decrease}{
  Let $\phi(s) = f_0 + \pd{g,s} + \half \pd{Hs,s}$ be a quadratic
  polynomial in $s\in \Re^n$, and $m(s) = \phi(s)+\sixth \sigma\nr{s}^3$, where $\sigma>0$ is a
  constant and $\nr{\cdot}$ is a general norm. Consider $s\neq 0$ and $u$ an
  eigenvector corresponding to $\lambda_{\min}[H]$ with $\|u\|_2=1$ and whose
  sign is chosen to ensure that $\pd{g+Hs,u}\leq 0$.
  Also assume that $\lambda_{\min}[H]+\sigma\nr{s} <0$.
  Then there exists an $\alpha >0$ such that
  \beqn{mdecr-2nd-2}
  m(s) - m(s+\alpha \nr{s} u)
  \geq \frac{3(\lambda_{\min}[H]+\sigma\nr{s})}{4\sigma^2}
       \left[\psi(s) \sigma^2\nr{s}^2 - \frac{3}{4}(\lambda_{\min}[H]+\sigma\nr{s})^2 \right],
  \eeqn
  where
  \beqn{psi-def}
  \psi(s) \eqdef \max\left[0,1+2\frac{\pd{g+Hs,u}}{\sigma\nr{s}^2}\right].
  \eeqn
}

\proof{
Setting $d = \nr{s}u$, we have that, for $\alpha > 0$,
\begin{align*}
m_k(s+ &\alpha d)
 =  m_k(s) + \alpha \pd{g+Hs,d} + \half \alpha^2\pd{Hd,d}
   + \sixth \sigma \nr{s+ \alpha d}^3-\sixth \sigma\nr{s}^3\\
& \le  m_k(s) + \half\alpha\sigma\nr{s}^3\,\left(2\frac{\pd{g+Hs,u}}{\sigma\nr{s}^2}\right)
   + \half \alpha^2\lambda_{\min}[H] \nr{s}^2
   + \sixth \sigma \nr{s+ \alpha s}^3-\sixth \sigma\nr{s}^3,
\end{align*}
where we have used the fact that $\nr{d} =\nr{s}$ implies the inequality
$\nr{s+ \alpha d}^3\le\nr{s+ \alpha s}^3$. Moreover 
\[
\nr{s+ \alpha s}^3-\nr{s}^3
= \big[(1+\alpha)^3-1\big]\,\nr{s}^3
=\alpha \big[ 3 + 3\alpha + \alpha^2 \big] \,\nr{s}^3,
\]
and hence, using \req{psi-def},
\beqn{mkd-noc}
m_k(s+ \alpha d)
\le  m_k(s) + \sixth \sigma\nr{s}^3\big(
  3\alpha\psi(s)
  + 3 \alpha^2 + \alpha^3 \big) + \half \alpha^2 \lambda_{\min}[H]\nr{s}^2.
\eeqn
This in turn yields that, for $\alpha > 0$,
\[
m(s)-m(s+\alpha d)
\geq -\frac{\alpha \nr{s}^2}{2}
     \left[\frac{\sigma\nr{s}}{3}\alpha^2 + (\lambda_{\min}[H]+\sigma\nr{s})\alpha +
       \sigma\nr{s}\psi(s)\right]
\eqdef  -\frac{\alpha \nr{s}^2}{2} q_0(\alpha).
\]
Now $q_0= a\alpha^2 + b\alpha +c$ is a convex quadratic in $\alpha$ which
admits a minimum for $\alpha = -b/(2a)$ of value $q(-b/(2a)) = c- b^2/(4a)$.
Since $b= \lambda_{\min}[H]+\sigma\nr{s} <0$, 
\[
m(s)-m(s+\alpha d)
\geq
\left(\frac{3(\lambda_{\min}[H]+\sigma\nr{s})}{2\sigma\nr{s}}\right)\frac{\nr{s}^2}{2}
\left[\sigma\nr{s}\psi(s) - \frac{3(\lambda_{\min}[H]+\sigma\nr{s})^2}{4\sigma\nr{s}}\right],
\]
for $\alpha = -b/(2a) >0$, which implies \req{mdecr-2nd-2}.
}

\noindent
This leads to the following necessary optimality condition.
  
\lthm{noc-2nd-order}{
  Let $\phi(s) = f_0 + \pd{g,s} + \half \pd{Hs,s}$ be a quadratic
  polynomial in $s\in \Re^n$, and assume that $s_*$ is a global
  minimizer of $m(s) = \phi(s)+\sixth \sigma\nr{s}^3$, where $\sigma>0$ is a
  constant and $\nr{\cdot}$ is a general norm. Let $u$ be an eigenvector
  corresponding to $\lambda_{\min}[H]$ such that $\|u\|_2=1$ and whose sign
  is chosen to ensure that $\pd{g+Hs_*,u}\leq 0$. 
  Then 
  \beqn{noc-2nd}
  \lambda_{\min}[H] + \omega(s_*) \sigma \nr{s_*} \ge 0,
  \eeqn
   where 
  \beqn{omega-def}
  \omega(s) \eqdef
  \left\{
  \begin{array}{ll}
  1+\bigfrac{2\sqrt{\psi(s)}}{\sqrt{3}}
  \le 1 + \bigfrac{2}{\sqrt{3}} \eqdef \kappa_\omega
  & \tim{if} s \neq 0,\\
  1 & \tim{otherwise,}
  \end{array}\right.
  \eeqn
  and $\psi(s)$ is given by \req{psi-def}.
}

\proof{When $H$ is positive-semidefinite, \req{noc-2nd} follows trivially. Assume
now that $H$ admits at least one negative eigenvalue.  Suppose first
that $s_*\neq 0$. If $\lambda_{\min}[H]+\sigma\nr{s_*}\geq 0$, \req{noc-2nd}
trivially follows.  Suppose thus that $\lambda_{\min}[H]+\sigma\nr{s_*} < 0$.
Then \req{mdecr-2nd-2} implies that there exists an $\alpha>0$ such that
$m(s_*+\alpha \nr{s_*} u) < m(s_*)$ (which is impossible), unless
\[
\psi(s_*) \sigma^2\nr{s_*}^2 > \frac{3}{4}(\lambda_{\min}[H]+\sigma\nr{s_*})^2.
\]
If $\psi(s_*)=0$, this cannot happen. Otherwise, this last inequality
requires that
\[
\sqrt{\psi(s_*)} \sigma\nr{s_*}
> \sqrt{\frac{3}{4}}\,\Big|\lambda_{\min}[H]+\sigma\nr{s_*}\Big|
> -\frac{\sqrt{3}}{2}\Big(\lambda_{\min}[H]+\sigma\nr{s_*}\Big),
\]
which, given \req{omega-def}, yields \req{noc-2nd}.

Suppose now that $s_*=0$ and that $\lambda_{\min}[H] <0$. It is then easy to verify
that, if the sign of $u$ is chosen to ensure that $\pd{g,u} \le 0$ and
\[
\alpha \in \left[0, -\frac{3\lambda_{\min}[H]}{2\sigma}\right],
\]
then
\beqn{DmE2}
m(\alpha u)
= f_0 + \alpha \pd{g,u} +\half \alpha^2\lambda_{\min}[H] + \sixth\alpha^3 \sigma
\le f_0+\frac{9\lambda_{\min}[H]^3}{16\sigma^2}
= m(0)+\frac{9\lambda_{\min}[H]^3}{16\sigma^2}
<m(0),
\eeqn
which again contradicts the assumption that $s_*=0$ is a
global minimum of $m$.  Thus $\lambda_{\min}[H] \ge 0$ and
\req{noc-2nd} also holds.
} 

\noindent
It is remarkable that this lemma provides a ``second-order'' necessary condition
for a global minimizer of quadratic polynomial regularized with a cubic term
in a possibly non-smooth norm, despite the first and second derivatives of
this objective function failing to exist.

\begin{figure*} 
\begin{center}
\centerline{
\includegraphics[width=6.8cm]{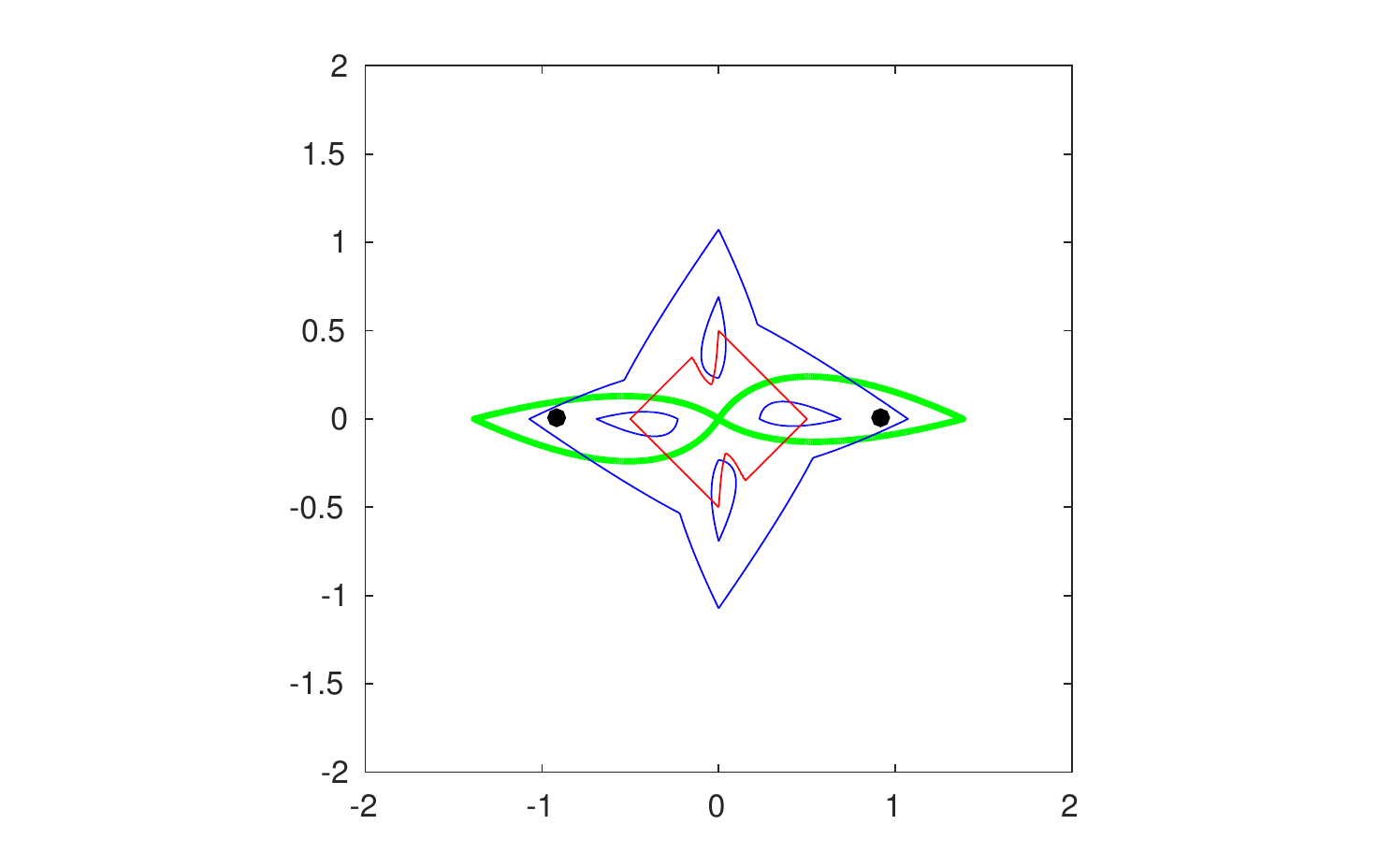}   
\hspace*{-1.4cm}
\includegraphics[width=6.8cm]{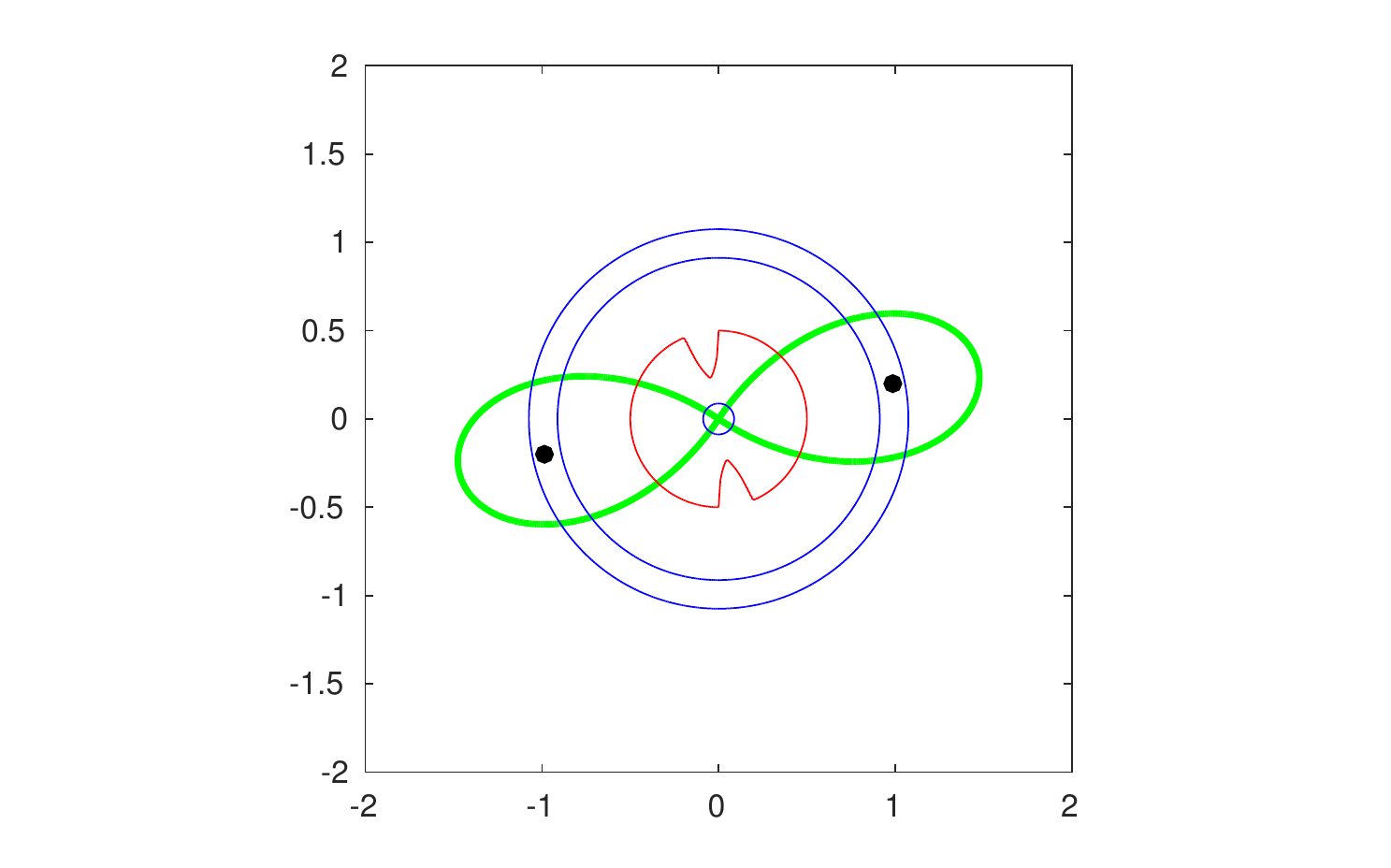}   
\hspace*{-1.4cm}
\includegraphics[width=6.8cm]{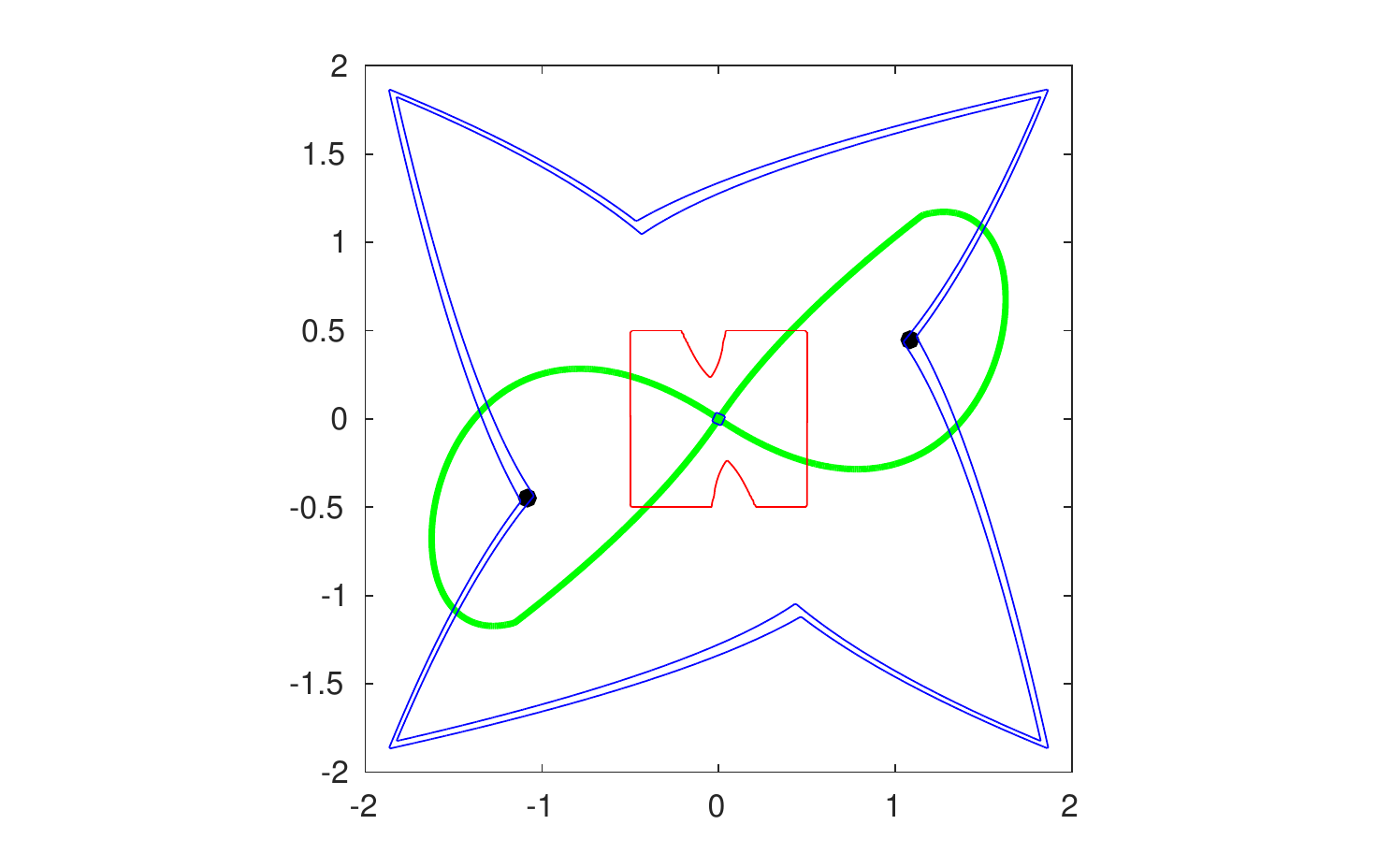} 
}
\caption{Admissible regions for the $\ell_1$- (left), $\ell_2$- (middle) and
  $\ell_\infty$- (right) norms} \label{fig:ell}
\end{center}
\vspace*{-1cm}
\end{figure*}

It is interesting to pause at this point to stress that the necessary first-
and second-order conditions \req{noc-1st} and \req{noc-2nd}, while sufficient
for our purposes as we will see, are \emph{merely necessary}, and by no means
sufficient to guarantee a local minimizer.  This is illustrated in
Figure~\ref{fig:ell}.  In this figure, a two-dimensional
model is constructed with a zero gradient at the origin and an indefinite
Hessian\footnote{Chosen to be $I-2uu^T/\pd{u,u}$ with $u^T= (5, 1)$.},
and with the regularization parameter $\sigma$ is
chosen equal to 6.  The left picture corresponds to the choice
$\nr{\cdot}=\|\cdot\|_1$, the central one to $\nr{\cdot}=\|\cdot\|_2$ and the
right one to $\nr{\cdot}=\|\cdot\|_\infty$, all other parameters being
identical. In each case, the region of the plane where \req{descent} holds is the interior
of the two green lobes and the regions were the deviation from \req{noc-1st} is
bounded by $0.01\,(\half)\sigma\nr{s}^2$ are shown in blue, the
first\footnote{In the central and left pictures.} being the small region
surrounding the origin (where the gradient is zero) and the second the zone
between the two blue concentric curves. Finally, the region where the
deviation from \req{noc-2nd} does not exceed $0.1\, \sigma\nr{s}$ is the
exterior of the region around the origin delineated in red. Note that, when it
exists, the region around the origin which is admissible for \req{noc-1st}
alone is excluded for \req{noc-2nd}. Thus in all cases, the admissible regions
for \req{descent}, \req{noc-1st} and \req{noc-2nd} consist of the regions
limited by any of the shown curves and containing the minimizers marked with a
black dot. We immediately notice that these regions are relatively large and
may extend reasonably far from the minimizers.  We also see that the geometry
of these regions, while simple for the Euclidean norm, can be quite
complicated for other norms.

Our algorithm for finding second-order $\epsilon_2$-approximate minimizers is
described \vpageref{AR2GN}.

\algo{AR2GN}{Second-Order Adaptive Regularization with General Norm (\tal{AR2GN})}{
  \begin{description}
  \item[Step 0: Initialization. ] An initial point $x_0\in \Re^n$, a regularization
    parameter $\sigma_0$ a desired final gradient accuracy
    $\epsilon, \epsilon_2 \in (0,1]$ and a model degree $p=2$ are given. The constants
    $\eta_1$, $\eta_2$, $\gamma_1$, $\gamma_2$, $\gamma_3$, $\theta_1>1$,
    $\theta_2>1$, and $\sigma_{\min}$
    are also given such that
    \beqn{eta-gamma2}
    \sigma_{\min} \in (0, \sigma_0], \;\;
      0 < \eta_1 \leq \eta_2 < 1, \tim{and}
      0< \gamma_1 < 1 < \gamma_2 < \gamma_3.
    \eeqn
    Compute $f(x_0)$ and set $k=0$.
  \item[Step 1: Check for termination. ]
    Compute $\lambda_{\min}[\nabla_x^2f(x_k)]$, and terminate with $x_\epsilon = x_k$ if 
    \beqn{term2}
    \|\nabla_x^1f(x_k)\|_{r,1} \leq \epsilon_1
    \tim{ and }
    \lambda_{\min}[\nabla_x^2f(x_k)] \ge - \epsilon_2.
    \eeqn
  \item[Step 2: Step calculation. ] Compute a step $s_k$  which
    sufficiently reduces the model $m_k$ in the sense that
    \req{descent} and \req{term-model} hold
    (for $p=2$)
    and, additionally,
    \beqn{term-model2}
    \lambda_{\min}[\nabla_x^2f(x_k)] + \theta_2\omega(s_k)\sigma_k\nr{s_k}\ge 0.
    \eeqn
  \item[Step 3: Acceptance of the trial point. ]
    Compute $f(x_k+s_k)$ and define $\rho_k$ as in \req{rhok-def}.
    If $\rho_k \geq \eta_1$, then define
    $x_{k+1} = x_k + s_k$; otherwise define $x_{k+1} = x_k$.
  \item[Step 4: Regularization parameter update. ]
   Set $\sigma_{k+1}$ according to \req{sigma-update}, increment $k$ by one and go to Step~1.
  \end{description}
}

\noindent
As is the case for $\theta_1$ in the \al{AR$1p$GN} algorithm, choosing an
iteration dependent $\theta_{2,k}$ is possible provided it is strictly bounded below by
$1$ and bounded above by a constant.

The existence of a suitable step in the \al{AR2GN} algorithm directly hinges
on Theorem~\ref{noc-2nd-order}.

\llem{step-exists2}{A step satisfying \req{descent}, \req{term-model}
  (for $p=2$)
  and \req{term-model2} always exists.
}

\proof{
It follows from Corollary~\ref{step-exists} and Theorem~\ref{noc-2nd-order} (with
$\phi(s) = T_{f,2}(x_k,s)$, $m=m_k$ and $\sigma = \sigma_k$) and the bound
$\theta_2 > 1$ that the required conditions are satisfied at
every global minimizer of the model $m_k$.
}

\noindent
As for the first-order case, continuity of $\nr{s}$ and of $T_{f,2}(x,s)$ with
respect to $s$ implies that conditions \req{descent}, \req{term-model}
and \req{term-model2} also hold in a neighbourhood
of a global minimizer whenever $\theta_1 > 1$ and $\theta_2>1$.
Such a neigbourhood can be reached for instance by using the algorithm
discussed in Section~\ref{RQMIN-s}.

Noting that Lemmas~\ref{SvsU}, \ref{lipschitz}, \ref{model-decrease} and
\ref{sigma-max} remain valid for the \al{AR2GN} algorithm, we
now provide a lower bound on the length of the step, which simplifies
that of \cite[Lemma~3.4]{CartGoulToin20}.

\llem{useful2}{
  Suppose that AS.1 holds for $p=2$ and that $k\in\calS$ before termination
  and such that $\lambda_{\min}[\nabla_x^2f(x_{k+1})]<-\epsilon_2$.  Then
  \[
  \nr{s_k} \geq  \frac{1}{L_{2,p}+ \theta_2\kappa_\omega\sigma_{\max}}\,\epsilon_2.
  \]
}

\proof{
Let $k\in \calS$ such that $\lambda_{\min}[\nabla_x^2f(x_{k+1})]<-\epsilon_2$. Since
$\min_z[a(z)+b(z)] \ge \min_za(z) + \min_zb(z)$, we deduce that
\begin{align*}
\lambda_{\min}[\nabla_x^2f(x_{k+1})] &
 =  \min_{\|d\|_2=1}\pd{\nabla_x^2f(x_{k+1})d,d} \\
& = \min_{\|d\|_2=1}\left[\pd{\nabla_x^2f(x_{k+1})d,d} -
  \pd{\nabla_x^2f(x_k)d,d}+\pd{\nabla_x^2f(x_k)d,d}\right]\\
&\ge \min_{\|d\|_2=1}\left[\pd{\nabla_x^2f(x_{k+1})d,d} -\pd{\nabla_x^2f(x_k)d,d}\right] + \min_{\|d\|_2=1}\pd{\nabla_x^2f(x_k)d,d}\\
& = \min_{\|d\|_2=1}\pd{(\nabla_x^2f(x_{k+1}) -\nabla_x^2f(x_k))d,d} + \min_{\|d\|_2=1}\pd{\nabla_x^2f(x_k)d,d}\\
&\ge - \|\nabla_x^2f(x_{k+1})-\nabla_x^2f(x_k)\|_2 + \lambda_{\min}[\nabla_x^2f(x_k)]\\
&\ge - \left(L_{2,p}\nr{s_k} + \theta_2 \omega(s_k)\sigma_k\nr{s_k}\right),
\end{align*}
where we also used \req{Lip-H} with $p=2$,
\req{term-model2} and \req{omega-def}. The
conclusion of the lemma then follows from Lemma~\ref{sigma-max} and the fact
that $\lambda_{\min}[\nabla_x^2f(x_{k+1})] < -\epsilon_2$.
} 

\noindent
We conclude our analysis by stating our final evaluation complexity bound for
finding second-order $\epsilon_2$-approximate minimizers.

\lthm{complexity2}{
  Suppose that AS.1--AS.2 hold
  for $p=2$
  and let
   \[
  \kappa_{\sf AR2GN} \eqdef
  \max\left\{
  \big[\half(L_{r,p}+\theta_1\sigma_{\max})\big]^{3/2}, 
      \left[L_{2,p}+\theta_2\kappa_{\omega}\sigma_{\max}\right]^{3}, 
      \right\}.
    \]
  Then the
  \al{AR2GN}
  algorithm requires at most
  \beqn{ar2gn-d}
  \left(\frac{6\,\kappa_{\sf AR2GN}}{\eta_1\sigma_{\min}}\right)
  \frac{f(x_0)-f_{\rm low}}{\min\left[\epsilon_1^{3/2},\epsilon_2^{3}\right]}
  \eeqn
  successful iterations and evaluations of $g$ and $H$  and at most
  \beqn{ar2gn-f}
  \left(\frac{6\,\kappa_{\sf AR2GN}}{\eta_1\sigma_{\min}}\right)
    \frac{f(x_0)-f_{\rm low}}{\min\left[\epsilon_1^{3/2},\epsilon_2^3\right]}
  \left(1+\frac{|\log\gamma_1|}{\log\gamma_2}\right)+
  \frac{1}{\log\gamma_2}\log\left(\frac{\sigma_{\max}}{\sigma_0}\right)
  \eeqn
  evaluations of $f$ to produce a vector $x_{\epsilon}\in \Re^n$ such that
  \[
  \|\nabla_x^1f(x_{\epsilon})\|_{r,1}\leq \epsilon_1 \tim{ and }
  \lambda_{\min}[\nabla_x^2f(x_{\epsilon})] \ge -\epsilon_2.
  \]
}

\proof{
  We prove the upper bounds \req{ar2gn-d} and \req{ar2gn-f} on the number of evaluations
  requested to produce an iterate $x_\epsilon$ at iteration $k_\epsilon$ such that 
  \[
  \|\nabla_x^1f(x_{\epsilon})\|_{r,1}\leq \epsilon_1 \tim{ and }
  \lambda_{\min}[\nabla_x^2f(x_\epsilon)] \ge -\epsilon_2
  \]
  is identical to that of Theorem~\ref{complexity}, except that $p=2$
  and the decrease
  \[
  \left(\frac{p!}{L_{r,p}+\theta_1\sigma_{\max}}\right)^{\sfrac{p+1}{p}}\epsilon_1^{\sfrac{p+1}{p}}
  \]
  is now replaced, using Lemmas~\ref{useful} and \ref{useful2}, by
  \[
  \min\left\{
  \left(\frac{2}{L_{r,p}+\theta_1\sigma_{\max}}\right)^{\sfrac{3}{2}} \epsilon_1^{\sfrac{3}{2}},
  \left(\frac{1}{L_{2,p}+\theta_2\kappa_\omega\sigma_{\max}}\right)^3 \epsilon_2^3
  \right\}.
  \]
  yielding the constant $\kappa_{\sf AR2GN}$. 
} 

\numsection{An algorithm for approximate minimization of regularized quadratics}
\label{RQMIN-s}

The rest of the paper is devoted to the definition and analysis of a method
whose purpose is to minimize a model of the form
\beqn{rq-def}
m(s) = f_0 + \pd{g,s}+\half \pd{Hs,s} + \sixth \sigma \nr{s}^3,
\eeqn
approximately, but enough for the conditions requested in Step~2 of the
\al{AR$1p$GN} (for $p=2$) and \al{AR2GN} algorithms to hold for $m=m_k$.
We first state a simple technical lemma.

\llem{root-bound}{
  Consider the quadratic polynomial $q(t)=at^2 + bt + c$ with $a\neq0$ and $c>0$.
  Then, for any $\nu>0$,
  \beqn{tstar}
  q(t_*) > \half c
  \tim{ for }
  t_* =  \min\left[ \,\frac{c}{\nu+3|b|},\,\frac{1}{3}\sqrt{\frac{c}{|a|}}\,\right].
  \eeqn
}

\proof{
  We immediately deduce that
  \[
  q(t_*) \ge c - |b|\left(\frac{c}{\nu+3|b|}\right) -|a|
  \left(\frac{1}{9}\frac{c}{|a|}\right)
  \ge c\Big( 1 - \sfrac{1}{3}-\sfrac{1}{9}\Big)
  > \half c.
  \]
} 

\noindent
The constant $\nu$ in \req{tstar}
is introduced to safeguard against $b=0$ and its value can be chosen for
convenience in what follows.

We may start building our specialized method (which we will call the
\al{RQMIN} algorithm) for minimizing the regularized quadratic \req{rq-def}.
The algorithm will unsurprisingly be iterative and we will denote its
successive iterates by $\{s_k\}_{k \ge 0}$ (the index $k$ refers, for the rest
of this section, to \al{RQMIN} iterations, and thus $g_k=g+Hs_k$). We will also make the choice to
start from the origin, that is $s_0=0$.  Moreover, we will construct the
iterates $s_k$ to ensure that the sequence $\{m(s_k)\}_{k\ge0}$ is
monotonically decreasing from $m(0)$. To motivate the forthcoming detailed
description of the algorithm ensuring this property, we now consider the
magnitude of the model decrease which can be obtained at a given iterate
$s_k$, if any.  We know from Lemma~\ref{noc-1st-order} that, if $s_k$ were a
local minimizer of $m$, then
\beqn{isequal}
\|g_k\|_{r,1} = \half\nr{s_k}^2.
\eeqn
This is the condition that the \al{RQMIN} algorithm will strive to achieve.
If \req{isequal} fails, we will now show that taking a step from
$s_k$ along a well chosen direction $d_k$ does produce a model decrease
\beqn{decr-1}
\Delta m(\alpha)
\eqdef m(s_k)-m(s_k+\alpha d_k)
=-\alpha \pd{g_k,d_k} - \half \alpha^2\pd{Hd_k,d_k}
- \sixth \sigma\nr{s_k+\alpha d_k}^3 + \sixth \sigma\nr{s_k}^3
\eeqn
which is suitably large. We start by analyzing the case
where the step is too short (in view of \req{isequal}), in which case a
generalized ``Cauchy point'' will provide adequate descent.

\llem{decrC}{
Let $s_k\in\Re^n$ such that $m(s_k)\le m(0)$ and
\beqn{ssmall}
\|g_k\|_{r,1}\geq \half\sigma\nr{s_k}^2.
\eeqn
Then
\beqn{DmC}
m(s_k) - m(s_k^C)
\geq \half \min\left[
  \frac{\big|\|g_k\|_{r,1}-\half\sigma\nr{s_k}^2\big|^2}{1+\sfrac{3}{2}\left(\|H\|_{r,2}+\sigma\nr{s_k}\right)},
  \frac{\big|\|g_k\|_{r,1}-\half\sigma\nr{s_k}^2\big|^{\sfrac{3}{2}}}{3\sqrt{\sigma}}
  \right],
\eeqn
where $g_k= g+Hs_k$ and
\beqn{mm-sC}
    s_k^C = s_k+\alpha_k^C d_k,
    \tim{with}
    d_k = \argmin_{\nr{v}=1}\pd{g_k,v}  \tim{and}
    \alpha_k^C = \argmin_{\alpha >0} m(s_k+\alpha d_k ).
\eeqn
}

\proof{If $\|g_k\|_{r,1}= \half\sigma\nr{s_k}^2$, the definition of $s_k^C$
  implies that $m(s_k)-m(s_k^C)\geq 0$ and \req{DmC} trivially follows.  Suppose
  therefore that the inequality in \req{ssmall} is strict, and consider the
  unidimensional minimization of $m(s_k+\alpha d_k)$ as a function of the scalar
  $\alpha$, giving \req{decr-1}.

Suppose first that $s_k=0$ and thus that $g_k=g$. Then
\[
\Delta m(\alpha) = \alpha\,q_0(\alpha)
\tim{ where }
q_0(\alpha) = \|g_k\|_{r,1}-\half \alpha \pd{Hd_k,d_k} - \sixth\sigma\alpha^2.
\]
We have that $q_0(0) =  \|g_k\|_{r,1}>0$ and $q_0(\alpha)$ is a concave
quadratic.  Hence the equation $q_0(\alpha)=0$ has a positive real root and we
may apply Lemma~\ref{root-bound} with $\nu= 1$ to deduce that
\[
q_0(\alpha_0) > \half\|g_k\|_{r,1}
\tim{ where }
\alpha_0 = \min\left[\frac{\|g_k\|_{r,1}}{1+\sfrac{3}{2}\left|\pd{Hd_k,d_k}\right|},
                    \frac{1}{3}\sqrt{\frac{\|g_k\|_{r,1}}{\sixth\sigma}}\right],
\]
and thus that
\beqn{case0}
\Delta m(\alpha_0)
\ge \alpha_0\,q_0(\alpha_0)
> \half \|g_k\|_{r,1} \min\left[\frac{\|g_k\|_{r,1}}{1+\sfrac{3}{2}\|H\|_{r,2}},
                    \frac{1}{3}\sqrt{\frac{\|g_k\|_{r,1}}{\sigma}}\right].
\eeqn

Suppose now that $\nr{s_k}>0$ and define
$v_k= \nr{s_k}d_k$. Then, because $\nr{v_k} = \nr{s_k}$, we have that 
$\nr{s_k+\alpha s_k} \ge \nr{s_k+\alpha v_k}$ and hence, from \req{decr-1},
\beqn{decr-2}
\Delta m(\alpha)
\ge \alpha \|g_k\|_{r,1}\nr{s_k} - \half \alpha^2\pd{Hv_k,v_k}
    - \sixth \sigma\nr{s_k+\alpha s_k}^3 + \sixth \sigma\nr{s_k}^3.
\eeqn
Observe now that
\[
\nr{s_k+\alpha s_k}^3 -\nr{s_k}^3
= [(1+\alpha)^3 - 1] \nr{s_k}^3
= \alpha(3 + 3\alpha + \alpha^2)\nr{s_k}^3,
\]
and thus \req{decr-2} becomes
\beqn{decr-3}
\Delta m(\alpha)
\ge \alpha \left[\|g_k\|_{r,1}\nr{s_k} - \half \alpha\pd{Hv_k,v_k}
  - \sixth \sigma(3+3\alpha+\alpha^2)\nr{s_k}^3\right]
\eqdef \alpha \, q_1(\alpha),
\eeqn
where
\[
q_1(\alpha) = \left(\|g_k\|_{r,1}\nr{s_k}-\beta\right)
- \alpha \left(\half \pd{Hv_k,v_k}+\beta\right)
-\alpha^2 \left( \third \beta \right),
\]
where we have defined $\beta \eqdef \half \sigma \nr{s_k}^3$.
Note that the constant term is positive because we have assumed
that $\|g_k\|_{r,1}>\half\sigma\nr{s_k}^2$. Moreover, $q_1(\alpha)$ is concave.
As above, this implies that the equation $q_1(\alpha)= 0$ has a positive real
root and we may then apply Lemma~\ref{root-bound} with $\nu = \half\nr{s_k}^2$
to deduce that
\[
q_1(\alpha_1) >\frac{1}{2}\left(\|g_k\|_{r,1}\nr{s_k}-\beta\right)
\]
where
\[
\alpha_1 \eqdef
\min\left[ \frac{\|g_k\|_{r,1}\nr{s_k}-\beta}{\half\nr{s_k}^2+3\left|\half\pd{Hv_k,v_k}+\beta\right|},
\frac{\sqrt{\|g_k\|_{r,1}\nr{s_k}-\beta}}{3\sqrt{\sfrac{1}{3}\beta}}\right],
\]
and hence, from \req{decr-3},
\begin{align*}
  \Delta m(\alpha_1)&
  \ge \alpha_1\,q(\alpha_1)\\
&\geq \frac{1}{2}\left(\|g_k\|_{r,1}\nr{s_k}-\beta\right)\min\left[
  \frac{\|g_k\|_{r,1}\nr{s_k}-\beta}{\half\nr{s_k}^2+3\left|\half\pd{Hv_k,v_k}+\beta\right|},
  \frac{\sqrt{\|g_k\|_{r,1}\nr{s_k}-\beta}}{3\sqrt{\third\beta}}\right]\\
&= \frac{1}{2}\min\left[ \frac{(\|g_k\|_{r,1}\nr{s_k}-\beta)^2}{\half\nr{s_k}^2+3|\half\pd{Hv_k,v_k}+\beta|},
  \frac{(\|g_k\|_{r,1}\nr{s_k}-\beta)^{\sfrac{3}{2}}}{3\sqrt{\third\beta}}\right].
\end{align*}
Using the identity $\nr{v_k} = \nr{s_k}$ and substituting the
definition of $\beta$, this finally gives that
\beqn{case1}
\Delta m(\alpha_1)
\ge \frac{1}{2}\min\left[ \frac{(\|g_k\|_{r,1}-\half\sigma\nr{s_k}^2)^2}{1+\sfrac{3}{2}(\|H\|_{r,2}+\sigma\nr{s_k})},
  \frac{(\|g_k\|_{r,1}-\half\sigma\nr{s_k}^2)^{\sfrac{3}{2}}}{3\sqrt{\sigma}}\right].
\eeqn

Combining now \req{case0} and \req{case1} gives \req{DmC}.
} 

\noindent
As is standard in Cauchy point approaches, the step $d_k$ in
\req{mm-sC} is made in the direction of the steepest descent for the
\emph{unregularized} quadratic, that is \emph{ignoring the regularization term}.

We now consider the alternative to \req{ssmall}, which, as \req{isequal}
indicates, means that the step $s_k$ is too large.  It therefore makes sense
to consider moving back from $s_k$ towards the origin.

\llem{decrR}{
Let $s_k\in\Re^n$ such that $m(s_k)\le m(0)$ and
\beqn{slarge}
\|g_k\|_{r,1} <  \half\sigma\nr{s_k}^2.
\eeqn
Then
\beqn{DmR}
m(s_k) - m(s_k^R)
\geq \half \min\left[
  \frac{\big|\|g_k\|_{r,1}-\half\sigma\nr{s_k}^2\big|^2}{1+\sfrac{3}{2}\left(\|H\|_{r,2}+\sigma\nr{s_k}\right)},
  \frac{\big|\|g_k\|_{r,1}-\half\sigma\nr{s_k}^2\big|^{\sfrac{3}{2}}}{3\sqrt{\sigma}}
  \right],
\eeqn
where $g_k= g+Hs_k$ and
\beqn{mm-sR}
    s_k^R = (1-\alpha_k^R) s_k
    \tim{with}
    \alpha_k^R = \argmin_{\alpha >0} m( s_k-\alpha s_k ).
\eeqn
}

\proof{
Note that \req{slarge} implies that $s_k\neq0$.
Then, from \req{decr-1} with $d_k = -s_k$,
\beqn{decr-2b}
\Delta m(\alpha)
\ge \alpha \pd{g_k,s_k} - \half \alpha^2\pd{Hs_k,s_k}
    - \sixth \sigma\nr{s_k-\alpha s_k}^3 + \sixth \sigma\nr{s_k}^3.
\eeqn
Now
\[
\nr{s_k-\alpha s_k}^3 -\nr{s_k}^3
= [(1-\alpha)^3 - 1] \nr{s_k}^3
= -\alpha(3 - 3\alpha + \alpha^2)\nr{s_k}^3,
\]
so that, from \req{decr-2b},
\beqn{decr-6}
\Delta m(\alpha)
\ge \alpha \left[\pd{g_k,s_k} - \half \alpha\pd{Hs_k,s_k}
  + \sixth \sigma(3-3\alpha+\alpha^2)\nr{s_k}^3\right]
\eqdef \alpha \, q_2(\alpha),
\eeqn
where
\[
q_2(\alpha) = \left(\pd{g_k,s_k}+\beta\right)
+ \alpha \left(-\half\pd{Hs_k,s_k}-\beta \right)
+ \alpha^2 \left(\third \beta \right).
\]
Observe now that, because  $s_k\neq 0$ and, since we have
assumed that $m(s_k) \le m(0)$, we have that $q_2(1)=\Delta m(1)
\le 0$. Moreover, the Cauchy-Schwarz inequality yields that 
\[
|\pd{g_k,s_k}|\le\|g_k\|_{r,1}\nr{s_k} < \half\sigma\nr{s_k}^3=\beta
\]
and hence 
\[
q_2(0) = \pd{g_k,s_k}+\beta > 0.
\]
This in turn implies the existence of a real root of $q_2(\alpha)$ in $(0,1]$,
and we may then again apply Lemma~\ref{root-bound}
with $\nu =\half\nr{s_k}^2$ to deduce that
\beqn{q2}
q_2(\alpha_2)
> \half \big(\pd{g_k,s_k}+\beta\big)
\eeqn
where
\beqn{alpha2}
\alpha_2
= \min\left[
  \frac{\pd{g_k,s_k}+\beta}{\half\nr{s_k}^2+3\left|-\half\pd{Hs_k,s_k}+\beta\right|},
  \frac{1}{3}\sqrt{\frac{\pd{g_k,s_k}+\beta}{\third\beta}}
  \right].
\eeqn
Moreover
\[
\pd{g_k,s}+\beta
\ge -\|g_k\|_{r,1}\nr{s} + \half \sigma\nr{s}^3
>0.
\]
Combining this bound with \req{decr-6}, \req{alpha2} and \req{q2}, we obtain that
\beqn{case2}
\Delta m(\alpha_2)
\geq \half \min\left[
  \frac{(\half\sigma\nr{s_k}^2-\|g_k\|_{r,1})^2}{1+\sfrac{3}{2}\left(\|H\|_{r,2}+\sigma\nr{s_k}\right)},
  \frac{(\half\sigma\nr{s_k}^2-\|g_k\|_{r,1})^{\sfrac{3}{2}}}{3\sqrt{\sigma}}
  \right],
\eeqn
which yields \req{DmR}.
} 

\noindent
Remarkably, \req{DmC} and \req{DmR} give identical lower bounds for the model decrease.
Lemmas~\ref{decrC} and \ref{decrR} generalize
\cite[Lemma~2.1]{CartGoulToin11} to the case where $s_k\neq0$ and general norms
are allowed.

We may now complete the analysis of what can happen at iterate $s_k$ (of the
still unspecified \al{RQMIN} algrithm) if the second-order necessary condition
of Theorem~\req{noc-2nd-order} fails. We first state an easy lemma giving
lower and upper bounds on the step $s_k$, dependent on the ``most negative
curvature'' of the quadratic given by \req{lar-def}.

\llem{slower}{Suppose that, for some $s_k$ and some $\beta\ge0$,
  \beqn{somedecr}
  m(s_0)-m(s_k) \ge \beta.
  \eeqn
  Then,
  \beqn{supper}
  \nr{s_k}
  \le \frac{\half \|H\|_{r,2}+
    \sqrt{\|H\|_{r,2}^2+\sfrac{2}{3}\sigma\|g\|_{r,1}}}{\third\sigma}
  \eqdef \kappa_{s,{\rm upp}},
  \eeqn
  and,  if $\beta>0$,
  \beqn{slow}
  \nr{s_k} \ge  \left\{\begin{array}{ll}
  \bigfrac{\sqrt{\|g\|_{r,1}^2+2\beta|\lambda_r[H]|}-\|g\|_{r,1}}{|\lambda_r[H]|}
    & \tim{if} \lambda_r[H] < 0\\
    \bigfrac{\beta}{\|g\|_{r,1}} & \tim{otherwise.}
    \end{array}\right.
  \eeqn
}

\proof{
  Since $s_0=0$ and $m(s_k) \ge m(0)+ \pd{g,s_k}+\half \pd{Hs_k,s_k}$, \req{somedecr}
  implies that
  \[
  -\|g\|_{r,1}\nr{s_k} + \half \min\big[0,\lambda_r[H]\big]\,\nr{s_k}^2
  \le \pd{g,s_k}+\half \pd{Hs_k,s_k}
  \le m(s_k) - m(0)
  \le -\beta,
  \]
  which gives \req{slow}.
  Observe now that \req{somedecr} implies that 
  \[
  \sixth \sigma\nr{s_k}^3
  \le |\pd{g,s_k}|+\half|\pd{Hs_k,s_k}|
  \le \|g\|_{r,1}\nr{s_k}+\half\|H\|_{r,2}\nr{s_k}^2,
  \]
 which yields \req{supper}.
} 

\noindent
Armed with this result, we now derive the model decrease when negative
curvature is present.

\llem{decrE}{
  Suppose that $u$ is an eigenvector associated with
  $\lambda_{\min}[H]$ with $\|u\|_2=1$ and that the sequence $\{m(s_k)\}_{k\ge0}$ is
  non-increasing.  For $k\ge 0$, define 
  \beqn{mm-sE}
  s_k^E = s_k + \alpha_k^E u_k
  \tim{where}
  u_k = - {\rm sign}\big(\pd{g_k,u}\big) u
  \tim{and}
  \alpha_k^E = \argmin_{\alpha >0} m(s_k+\alpha u_k ).
  \eeqn
  Then
  \beqn{DmE1b}
  m(s_0)-m(s_1) \ge m(s_0) - m(s_0^E) \geq \frac{9|\lambda_{\min}[H|^3}{16\sigma^2}
  \eeqn
  and there exists a constant $\kappa_s$ such that, for $k\ge 1$,
  \beqn{sbound}
  \nr{s_k} \ge \kappa_s.
  \eeqn
  Moreover, if
  \beqn{not-2nd}
    \lambda_{\min}[H] + \sigma\omega(s_k)\nr{s_k} <0
  \eeqn
  at iteration $k \ge 1$, then one has that
  \beqn{DmE2b}
  m(s_k) - m(s_k^E)
  \ge \bigfrac{9\nr{s_k}^2}{16\sigma^2}\,\Big|\lambda_{\min}[H] + \sigma\omega(s_k)\nr{s_k}\Big|^3.
  \eeqn
}

\proof{
The first inequality in \req{DmE1b} results from \req{mm-s} and the
second is a direct consequence of the proof of Theorem~\ref{noc-2nd-order}
(see \req{DmE2}). The existence of $\kappa_s$ such that \req{sbound} holds
for $k\ge 1$ then follows from Lemma~\ref{slower} with
$\beta= \sfrac{9}{16}|\lambda_{\min}[H]|^3/\sigma^2$
and our assumption that $\{m(s_k)\}$ is non-increasing.
We now prove \req{DmE2b}.  From \req{not-2nd}, we have that
\beqn{not-2nd-0}
\lambda_{\min}[H] + \sigma\omega(s_k)\nr{s_k} = -\mu \sigma\nr{s_k}
\eeqn
for some $\mu>0$. But
\req{omega-def} implies that
$0>\lambda_{\min}[H]+\sigma\nr{s_k} = -|\lambda_{\min}[H] +
\sigma\nr{s_k}|$, and thus, from \req{omega-def},
\beqn{not-2nd-2}
|\lambda_{\min}[H] + \sigma\nr{s_k}|
  = \frac{2\sqrt{\psi(s_k)}}{\sqrt{3}}\sigma\nr{s_k} + \mu \sigma \nr{s_k},
\eeqn
from which we obtain that
\[
\big(\lambda_{min}[H] + \sigma\nr{s_k}\big)^2
= \frac{4}{3}\psi(s_k)\sigma^2\nr{s_k}^2
  + \sigma^2\nr{s_k}^2\left(\mu^2 +
  \frac{4\sqrt{\psi(s_k)}}{\sqrt{3}}\,\mu\right).
\]
Substituting this inequality in \req{mdecr-2nd-2}, then gives that there
exists an $\alpha>0$ such that
\begin{align*}
m(s_k)-m(s_k +\alpha \nr{s_k} u)
&\geq \frac{3(\lambda_{\min}[H] + \sigma\nr{s_k})}{4\sigma^2}
\left[ -\frac{3}{4}\sigma^2\nr{s_k}^2
  \left(\mu^2 + \frac{4\sqrt{\psi(s_k)}}{\sqrt{3}}\,\mu\right)\right]\\
& = \frac{9}{16}\Big|\lambda_{\min}[H] + \sigma\nr{s_k}\Big|\nr{s_k}^2
\left(\mu^2 + \frac{4\sqrt{\psi(s_k)}}{\sqrt{3}}\,\mu\right).
\end{align*}
But \req{not-2nd-2} implies that $\big|\lambda_{\min}[H] + \sigma\nr{s_k}\big| \ge \mu
\sigma\nr{s_k}$, and thus
\[
m(s_k)-m(s_k +\alpha \nr{s_k} u)
\ge \frac{9\sigma\nr{s_k}^3}{16}
\left(\mu^3 + \frac{4\sqrt{\psi(s_k)}}{\sqrt{3}}\,\mu^2\right)
\ge \frac{9\sigma\nr{s_k}^3}{16} \mu^3.
\]
The inequality \req{DmE2b} then follows from \req{not-2nd-0} and \req{mm-sE}.
} 

\noindent
We now have all ingredients to describe the \al{RQMIN} algorithm.
It hinges on \req{decrC}, \req{decrR} and \req{decrE} and proceeds by
successive one-dimensional minimizations of $m$ along the directions $s_k^C$
or $s_k^R$ (depending on the sign of $\|g_k|_{r,1}- \half \sigma \nr{s}^2$)
and, if needed, $s_k^E$. It is formally stated \vpageref{mod-min}. 

\algo{mod-min}{An algorithm for minimization of a regularized quadratic (\tal{RQMIN})}{
  The value $f_0$, gradient $g$ and Hessian $H$ of the quadratic at $s=0$
  are given, as well as a regularization parameter $\sigma$ and accuracy
  requests $\epsilon_1>0$ and $\theta_2>1$.
  \begin{description}
  \item[Step 0: Initialization]
    If unavailable, compute $\lambda_{\min}[H]$ and and associated eignevector
    $u$ with $\|u\|_2=1$.
    Set $k = 0$, $s_0=0$ and $g_0=g$.
  \item[Step~1: Check for termination.]
    Terminate if
    \beqn{RQMIN-term}
    \left|\|g_k\|_{r,1} - \frac{\sigma}{2}\nr{s_k}^2\right| \le \epsilon_1
    \tim{ and }
    \lambda_{\min}[H] + \theta_2 \omega(s_k)\sigma \nr{s_k} \ge 0.
    \eeqn
  \item[Step~2: Negative gradient step.]
    If $\|g_k\|_{r,1} > \frac{\sigma}{2}\nr{s_k}^2$, compute $s_k^C$ according
    to \req{mm-sC}, set $m_{k,1}=m(s_k^C)$  and go to Step~4.
  \item[Step~3: Retraction step.]
    If $\|g_k\|_{r,1} < \frac{\sigma}{2}\nr{s_k}^2$, compute $s_k^C$ according
    to \req{mm-sR} and set $m_{k,1}=m(s_k^R)$. 
  \item[Step~4: Eigenvalue step.]
    If $\lambda_{\min}[H] +\theta_2\omega(s_k)\sigma\nr{s_k}<0$, compute
    $s_k^E$ according to \req{mm-sE} and set $m_{k,2} = m(s_k^E)$.
    Else, set $m_{k,2} = m(s_k)$.
  \item[Step~5: Move.]
    Set
    \beqn{mm-s}
    s_{k+1} = \left\{\begin{array}{ll} s_k^C  & \tim{if} m_{k,1} \le m_{k,2},\\
                                      s_k^E  & \tim{otherwise,}
                    \end{array}\right.
    \tim{and}
    g_{k+1} = g_k + H(s_{k+1}-s_k).
    \eeqn
    Increment $k$ by one and got to Step~1.
  \end{description}
} 

\noindent
Note that the mechanism of the algorithm, which proceeds by
successive unidimensional minimizations, guarantees that the sequence
$\{m(s_k)\}$ is monotonically decreasing, as announced. 

Having established, in Lemma~\ref{decrC}, \ref{decrR} and \ref{decrE}, lower
bounds on the decrease in $m$ for all steps produced 
by the \al{RQMIN} algorithm, we are now ready to state its iteration
complexity\footnote{At variance with its evaluation complexity, which would be
irrelevant here since evaluating $m(s)$ as many times as necessary does not
require evaluating $f_0$, $g$ and $H$ more than once (when the algorithm is called).}.

\lthm{rqmin-complexity}{Given $\epsilon_1 > 0$ and $\theta_2 > 1$,
  there exist a constant $\kappa_{\sf RQMIN}>0$ independent of $k$ such that
  the \al{RQMIN} algorithm requires at most 
  \beqn{sub-itbnd}
  \kappa_{\sf RQMIN} \max\left[\epsilon_1^{-2}, \epsilon_1^{-\sfrac{3}{2}},(\theta_2-1)^{-3}\right]
  \eeqn
  iterations to produce an iterate $s_k$ such that
  \beqn{term-conds}
  \left|\|g_k\|_{r,1} - \half\sigma\nr{s_k}^2\right| \le \epsilon_1
  \tim{and}
  \lambda_{\min}[H] + \theta_2 \omega(s_k)\sigma\nr{s_k} \ge 0.
  \eeqn
}

\proof{
 If the \al{RQMIN} algorithm terminates at $k=0$, then the bound
 \req{sub-itbnd} is trivially satisfied.  Assume therefore that termination
 does not occur at $s_0$. We therefore have that, for $k \ge 1$ before termination,
 \[
 \tim{either} \left|\|g_k\|_{r,1} - \half \theta_1\sigma\nr{s_k}^2\right| > \epsilon_1
 \tim{or}     \lambda_{\min}[H] +\theta_2 \omega(s_k) \sigma \nr{s_k}<0.
 \]
 Let us define $\calN \eqdef\{k\ge 0 \mid s_k = s_k^E\}$ and note that, by construction, this
 set is non-empty only if $\lambda_{\min}[H]<0$.
 We then obtain from $g_0=g$, \req{mm-s}, \req{DmC} and
 \req{DmE1b} that
 \beqn{dd}
 m(s_0)-m(s_1)
 \ge \left\{\begin{array}{ll}
  \half \max\left\{
  \min\left[\bigfrac{\|g\|_{r,1}^2}{1+\sfrac{3}{2}\|H\|_{r,2}},
  \bigfrac{\|g\|_{r,1}^{\sfrac{3}{2}}}{3\sqrt{\sigma}}\right],
  \bigfrac{9|\lambda_{\min}[H]|^3}{16\sigma^2}
  \right\}
  & \tim{ if } 0 \in \calN,\\*[3ex]
  \half \min\left[
  \bigfrac{\|g\|_{r,1}^2}{1+\sfrac{3}{2}\|H\|_{r,2}},
  \bigfrac{\|g\|_{r,1}^{\sfrac{3}{2}}}{3\sqrt{\sigma}},
  \right]
  & \tim{ otherwise.}\\
  \end{array}\right.   
 \eeqn
 Observe now that the second part of \req{term-conds} cannot hold as long as
 \[
 -\big(\lambda_{\min}[H]+ \omega(s_k)\sigma\nr{s_k}\big) >
 (\theta_2-1)\omega(s_*)\sigma\nr{s_k}.
 \]
 Hence \req{DmE2b} give that, for $k\ge 1$ and $k\in\calN$ before termination,
 \beqn{decr-sE}
 m(s_k)-m(s_k^E) \ge \frac{9\nr{s_k}}{16\sigma^2}\left[(\theta_2-1)\omega(s_*)\right]^3.
 \eeqn
 Because Lemma~\ref{slower} with $\beta$ chosen as the relevant right-hand side in
 \req{dd} guarantees the  existence of $\kappa_{s,{\rm low}}>0$ such that 
 $\nr{s_k} \ge \kappa_{s,{\rm low}}$ for all $k \ge 1$, and, because
 $\omega(s_k)\ge 1$, \req{decr-sE} ensures that, before termination and for
 $k\in\calN$, $k\ge 1$,
 \[
 m(s_k)-m(s_k^E) \ge \frac{9\kappa_{s,{\rm low}}^2}{16\sigma^2}(\theta_2-1)^3
 \,\eqdef \, \frac{9\kappa_{s,{\rm low}}^2}{16\sigma^2}\epsilon_2^3.
 \]
 Using this together with \req{mm-s} and \req{DmC} gives that
 for $k\ge 1$ before termination,
 \begin{align*}
 m(s_k)&-m(s_{k+1}) \ge \\
 &\left\{\begin{array}{ll}\half \max\left\{
  \min\left[\bigfrac{\epsilon_1^2}{1+\sfrac{3}{2}\left(\|H\|_{r,2}+\sigma\nr{s_k}\right)},
  \bigfrac{\epsilon_1^{\sfrac{3}{2}}}{3\sqrt{\sigma}}\right],
  \bigfrac{9\kappa_{s,{\rm low}}^2}{16\sigma^2}
  \right\}
  & \tim{if} k\in\calN,\\*[3ex]
  \min\left[\bigfrac{\epsilon_1^2}{1+\sfrac{3}{2}\left(\|H\|_{r,2}+\sigma\nr{s_k}\right)},
  \bigfrac{\epsilon_1^{\sfrac{3}{2}}}{3\sqrt{\sigma}}\right]
  & \tim{otherwise,}
 \end{array}\right.
 \end{align*}
  and therefore, using \req{supper},
  \begin{align}
  m(s_k)-m(s_{k+1})
  & \ge \half \min\left[
  \frac{1}{1+\sfrac{3}{2}\left(\|H\|_{r,2}+\sigma\kappa_{s,{\rm upp}}\right)},
  \frac{1}{3\sqrt{\sigma}},
  \bigfrac{9\kappa_{s,{\rm low}}^2}{16\sigma^2}
  \right]
  \,\min\Big[\epsilon_1^2,\epsilon_1^{\sfrac{3}{2}}, \epsilon_2^3\Big] \nonumber\\
  & \eqdef \kappa_* \,\min\Big[\epsilon_1^2,\epsilon_1^{\sfrac{3}{2}},\epsilon_2^3\Big]. \label{mod-decr}
  \end{align}
  We now observe that the definition of $m(s)$ and \req{supper} together
  imply that
  \[
  m(s)
  \ge m(0) - \|g\|_{r,1}\kappa_{s,{\rm upp}} - \half
  \|H\|_{r,2}\kappa_{s,{\rm upp}}^2
  \eqdef m_{\rm low}.
  \]
  Therefore \req{mod-decr} implies that the number of iterations required by
  the \al{RQMIN} algorithm to produce an iterate such that \req{term-conds}
  holds cannot exceed 
  \[
  \frac{m(0) - \beta-m_{\rm low}}
       {\kappa_*\min\Big[\epsilon_1^2,\epsilon_1^{\sfrac{3}{2}},\epsilon_2^3\Big]}
  \]
  which is \req{sub-itbnd} with $\kappa_{\sf RQMIN} = (m(0) - \beta-m_{\rm low})/\kappa_*$.
}

\noindent
We now consider applying  the \al{RQMIN} algorithm to find a step $s_k$ in
Step~2 of the \al{AR2GN} method\footnote{With $f(x_k)=f_0$, $\nabla_x^1f(x_k)= g$,
$\nabla_x^2 f(x_k)= H$, $\sigma_k = \sigma$ and $x_k=0$.}.
This latter methods requires the conditions \req{descent}, \req{term-model}
and \req{term-model2} to hold. We immediately note that \req{descent}
automatically holds because of the monotonically decreasing nature of the values of
$m$ in the \al{RQMIN} algorithm.  Moreover, \req{term-model2} and the second
part of \req{RQMIN-term} are identical.  However, the first part of
\req{RQMIN-term} is too strong, because it imposes a two-sided inequality
on $\|g_k\|_{r,1}-\half \sigma_k\nr{s_k}$ while \req{term-model} only requests
\beqn{weaker}
\|g_k\|_{r,1}-\half\sigma_k\nr{s_k}^2
\le \epsilon_{1s}
\eqdef \half(\theta_1-1)\sigma_k\nr{s_k}^2
\eeqn
but allows for $\|g_k\|_{r,1}-\half\sigma_k\nr{s_k}^2$ to be negative.
In Figure~\ref{fig:ell}, this amounts to removing the outer blue curve,
thus enlarging the admissible regions containing the minimizers.
A modified variant of the \al{RQMIN} algorithm is therefore suitable if our only
objective is to satisfy \req{descent}, \req{term-model} and \req{term-model2}.
This variant, which we call the \al{RQMIN1} algorithm, differs from
\al{RQMIN} in that
\begin{enumerate}
\item the first part of \req{RQMIN-term} is replaced by requiring that
  \req{weaker} holds,
\item Step~3 of \al{RQMIN} is skipped (as there is no need to correct for negative
      $\|g_k\|_{r,1}-\half\sigma_k\nr{s_k}^2$).
\end{enumerate}
In addition, because \req{weaker} is weaker that the first part of
\req{RQMIN-term}, termination of the \al{RQMIN1} algorithm cannot happen
later than that what would happen if applying the \al{RQMIN} algorithm with
$\epsilon_1=\epsilon_{1s}$. This allows us to derive the following upper bound
on the number of iterations of the \al{RQMIN1} algorithm that are necessary to
compute a step $s_k$ in Step~2 of \al{AR2GN}.

\lcor{RQMIN1-complexity}{
Given $\theta_1 > 1$ and $\theta_2 > 1$,
  there exist a constant $\kappa_{\sf RQMIN1}>0$ independent of $k$ such that the \al{RQMIN1} algorithm requires at most
  \beqn{sub-itbnd-1}
  \kappa_{\sf RQMIN1} \max\left[(\theta_1-1)^{-2},
                              (\theta_1-1)^{-\sfrac{3}{2}},
                              (\theta_2-1)^{-3}\right]
  \eeqn
  iterations to produce an iterate $s_k$ such that \req{descent},
  \req{term-model} and \req{term-model2} hold.
}

\proof{
  The desired result immediately follows by noting that, by virtue of
  \req{sigma-update}, \req{weaker} and the definition of $\kappa_{s,{\rm low}}$ in the
  proof of Theorem~\ref{rqmin-complexity},
  \[
  \epsilon_{1s} \geq \half (\theta_1-1)\sigma\kappa_{s,{\rm low}}^2.
  \]
  The bound \req{sub-itbnd-1} then follows with
  \[
  \kappa_{\sf RQMIN1}
  \eqdef \kappa_{\sf RQMIN}
    \min\left[(\half\sigma\kappa_{s,{\rm low}}^2)^2,
              (\half\sigma\kappa_{s,{\rm low}}^2)^{\sfrac{3}{2}}\right].
  \]
} 

\noindent
The reader may have wondered why we did consider the \al{RQMIN} method and
its two-sided condition at all, since its one-sided version \al{RQMIN1} is
sufficient for the purpose of computing a step in \al{AR2GN}.
Our motivation for \al{RQMIN} is that it is likely to achieve a larger model decrease, hopefully
reducing the number of iterations needed by \al{AR2GN} to terminate.

But the story does not finish here. As we have alluded to in Section~\ref{algo-s},
an even simpler variant of the \al{RQMIN} algorithm can be used to compute
$s_k$ in Step~2 of the \al{AR$1p$GN} algorithm when $p=2$. Since the only
requirements on $s_k$ are then \req{descent} and \req{term-model}, we may
define the \al{RQMIN2} algorithm as a variant of \al{RQMIN} where
\begin{enumerate}
\item the whole of \req{RQMIN-term} is replaced by requiring that
  \req{weaker} holds,
\item Step~3 and Step 4 of \al{RQMIN} are skipped and the first part
      of \req{mm-s} replaced by $s_{k+1} = s_k^C$.
\end{enumerate}
Removing all bounds related to the second-order condition in
Theorem~\ref{rqmin-complexity}, we then obtain the following iteration bound
for the \al{RQMIN2} algorithm (as needed in Step~2 of \al{AR$1p$GN} with $p=2$).

\lcor{RQMIN2-complexity}{
Given $\theta_1 > 1$,
  there exist a constant $\kappa_{\sf RQMIN2}>0$ independent of $k$ such that
  the \al{RQMIN2} algorithm requires at most 
  \beqn{sub-itbnd-2}
  \kappa_{\sf RQMIN2} \max\left[(\theta_1-1)^{-2},
                              (\theta_1-1)^{-\sfrac{3}{2}}\right]
  \eeqn
  iterations to produce an iterate $s_k$ such that \req{descent},
  and \req{term-model} hold.
}

\noindent
Note that the \al{RQMIN2} algorithm reduces to a standard first-order method
(in the $\nr{\cdot}$ norm), but applied to the quadratic alone, instead of to the
complete regularized model.

\numsection{Discussion}\label{conclusion-s}

We have shown in this variant of \cite{GratToin21} that using Euclidean
curvature in second-order necessary optimality conditions (see \req{lar-def}) leads to results
that are entirely similar to those obtained when attempting to enforce the
condition $\min_{\nr{v}=1}\pd{Hv,v} > \epsilon_2$ used in this reference.

{\footnotesize
  
\section*{{\footnotesize Acknowledgements}}

The authors are indebted to Sadok Jerad for his careful reading of the manuscript.
  

}

\appendix
\renewcommand{\theequation}{A.\arabic{equation}}
\renewcommand{\thesection}{A.\arabic{section}}
\renewcommand{\thesubsection}{A.\arabic{subsection}}

\end{document}